\documentclass[preprint,3p,11pt]{elsarticle}

\usepackage{tabls}
\usepackage{cites}
\usepackage{epsf}
\usepackage{appendix}
\usepackage{ragged2e}
\usepackage{enumitem}
\setlist[itemize]{leftmargin=*}
\usepackage[ruled,vlined]{algorithm2e}
\usepackage[caption = false]{subfig}
\usepackage{amsmath}
\usepackage{amssymb}
\usepackage{xcolor}
\usepackage[normalem]{ulem}

\usepackage{hyperref}
\hypersetup{colorlinks=true}

\newcommand{\dv}[2]{\frac{\partial #1}{\partial #2}}
\newcommand{\ip}[1]{\langle #1 \rangle}

\begin{document}
\begin{frontmatter}



\title{A Nonlinear Projection-Based Iteration Scheme with Cycles over
Multiple Time Steps for Solving Thermal Radiative Transfer Problems}

\author[lanl]{Joseph M. Coale}
\ead{jmcoale@lanl.gov}
\author[ncsu]{Dmitriy Y. Anistratov}
\ead{anistratov@ncsu.edu}
\address[lanl]{Los Alamos National Laboratory,
  Los Alamos, NM  87545}
\address[ncsu]{Department of Nuclear Engineering,
North Carolina State University,
 Raleigh, NC 27695}

\begin{abstract}
In this paper we present a multilevel projection-based iterative scheme for solving thermal radiative transfer problems that performs iteration cycles
on the high-order Boltzmann transport equation (BTE) and low-order moment equations.
Fully implicit temporal discretization based on the backward Euler time-integration method is used for all equations.
The multilevel iterative scheme is designed
to perform iteration cycles over collections of multiple time steps,
each of which can be interpreted as a coarse time interval  with a subgrid of time steps.
This treatment is demonstrated to transform implicit temporal integrators to diagonally-implicit multi-step schemes
on the coarse time grid formed with the amalgamated time intervals.
A multilevel set of moment equations are formulated by the
nonlinear projective approach.
The Eddington tensor defined with
the BTE solution provides exact closure for the moment equations.
During each iteration, a number of chronological time steps are solved with the BTE alone, after which the same collection of time steps is solved with the moment equations and material energy balance.
Numerical results are presented to demonstrate the effectiveness of this iterative scheme for simulating evolving radiation and heat waves in 2D geometry.
\end{abstract}


 \end{frontmatter}


\section{Introduction}
We consider the basic thermal radiative transfer (TRT) problem that neglects material motion, heat conduction and the scattering of photons which is formulated with the multigroup Boltzmann transport equation (BTE) describing photons
\begin{subequations}
\begin{equation}
		\frac{1}{c}\dv{I_g}{t} + \boldsymbol{\Omega}\cdot\boldsymbol{\nabla}I_g + \varkappa_g(T)I_g = \varkappa_g(T)B_g(T), \label{bte}
\end{equation}
\begin{equation}
		I_g|_{\boldsymbol{r}\in\partial\Gamma}=I_g^\text{in}\ \ \text{for}\ \ \boldsymbol{\Omega}\cdot\boldsymbol{n}_\Gamma<0,\quad I_g|_{t=0}=I_g^0,
\end{equation}
\label{eq:bte}
\end{subequations}
\begin{equation*}
		\boldsymbol{r}\in\Gamma,\quad
		t \in [0,t^\text{end}], \quad \boldsymbol{\Omega}\in\mathcal{S},\quad g=1,\dots,G \nonumber
\end{equation*}
\vspace{-0.1cm}
coupled to the material energy balance (MEB) equation that models
radiation-matter energy exchange
\begin{equation} \label{meb}
	\dv{\varepsilon(T)}{t} = { \sum_{g=1}^{G}\bigg(\int_{4\pi}I_gd\Omega - B_g(T)\bigg)\varkappa_g(T)},\quad T|_{t=0}=T^0,
\end{equation}
where $\boldsymbol{r}$ is spatial position, $\boldsymbol{\Omega}$ is the direction of particle motion, $g$ is the photon frequency group index, $c$ is the speed of light, $\Gamma$ is the spatial domain, $\partial\Gamma$ is the domain boundary, $\boldsymbol{n}_\Gamma$ is the unit normal to $\partial\Gamma$ and $\mathcal{S}$ is the unit sphere. $I_g(\boldsymbol{r},\boldsymbol{\Omega},t)$ is the group specific intensity of radiation, $T(\boldsymbol{r},t)$ is the material temperature, $\varkappa_g(\boldsymbol{r},t;T)$ is the group photon opacity, $\varepsilon(\boldsymbol{r},t;T)$ is the material energy density and $B_g(\boldsymbol{r},t;T)$ is the group Planck black-body distribution function.

The iterative algorithm is based on the multilevel quasidiffusion (MLQD) methodology
\cite{gol'din-1964,auer-mihalas-1970,PASE-1986,dya-aristova-vya-mm1996,aristova-vya-avk-m&c1999,dya-jcp-2019}, which is formulated via the nonlinear projective approach. Two systems of moment equations of the BTE are constructed by projecting the BTE onto a series of low-order subspaces.
The first system of moment equations is the multigroup low-order quasidiffusion (LOQD) equations
(aka Variable Eddington Factor equations), which are the first two angular moments of the BTE
\begin{subequations} \label{eq:mg-loqd}
	\begin{gather}
		\dv{E_g}{t} + \boldsymbol{\nabla}\cdot\boldsymbol{F}_g + c\varkappa_g(T)E_g = 4\pi\varkappa_g(T)B_g(T),\\
		\frac{1}{c}\dv{\boldsymbol{F}_g}{t} + c\boldsymbol{\nabla}\cdot (\boldsymbol{\mathfrak{f}}_g E_g) + \varkappa_g(T) \boldsymbol{F}_g = 0,
	\end{gather}
\end{subequations}
whose unknowns include the group radiation energy density $E_g=\frac{1}{c}\int_{4\pi}I_gd\Omega$ and flux $\boldsymbol{F}_g = \int_{4\pi}\boldsymbol{\Omega}I_gd\Omega$. The second system of moment equations is the effective grey LOQD equations obtained by summing Eqs. \eqref{eq:mg-loqd} over all frequency groups
\begin{subequations} \label{eq:gr-loqd}
	\begin{gather}
		\dv{E}{t} + \boldsymbol{\nabla}\cdot\boldsymbol{F} + c\ip{\varkappa}_E E = c\ip{\varkappa}_B a_RT^4,\\
		\frac{1}{c}\dv{\boldsymbol{F}}{t} + c\boldsymbol{\nabla} \cdot (\ip{\boldsymbol{\mathfrak{f}}}_{E}E) +\ip{\varkappa}_F\boldsymbol{F} + \bar{\boldsymbol{\eta}}E = 0,
	\end{gather}
\end{subequations}
which solve for the total radiation energy density $E=\sum_{g=1}^GE_g$ and flux $\boldsymbol{F}=\sum_{g=1}^G\boldsymbol{F}_g$.
The MEB equation \eqref{meb} is coupled with the effective grey LOQD equations \eqref{eq:gr-loqd} on the grey scale, reducing the dimensionality of the TRT problem. This is done by reformulating the MEB in effective grey form
\begin{equation} \label{grey_meb}
	\dv{\varepsilon(T)}{t} = c\ip{\varkappa}_E E - c\ip{\varkappa}_Ba_RT^4.
\end{equation}

The Eddington tensor $\boldsymbol{\mathfrak{f}}_g$  computed by the solution of the BTE   \eqref{eq:bte}
and
spectrum averaged coefficients calculated by the solution of multigroup LOQD equations \eqref{eq:mg-loqd}
which define exact closures for the LOQD system are
\begin{gather}
	\boldsymbol{\mathfrak{f}}_g = \frac{\int_{4\pi}\boldsymbol{\Omega}\otimes\boldsymbol{\Omega}{I}_g\ d\Omega}{\int_{4\pi}I_gd\Omega}, \label{eq:et}\\[2.5pt]
	\ip{u}_E = \frac{\sum_{g=1}^Gu_gE_g}{\sum_{g=1}^GE_g},\quad
	\ip{u}_F=\text{diag}(\ip{u_g}_{F_x},\ip{u_g}_{F_y},\ip{u_g}_{F_z}), \label{eq:grip2}\\[2.5pt]
	\ip{u}_{F_\alpha}=\frac{\sum_{g=1}^Gu_g|F_{\alpha,g}|}{\sum_{g=1}^G|F_{\alpha,g}|}, \quad
	\bar{\boldsymbol{\eta}}=\frac{\sum_{g=1}^{G} (\varkappa_g-\bar{\mathbf{K}}_{R})\boldsymbol{F}_g}{\sum_{g=1}^GE_g} \, .
	\label{eq:grip}
\end{gather}
\pagebreak

\noindent
The MLQD method for TRT problems is defined
by the  system of equations consisting of the following parts:
\begin{itemize}
\item the  high-order multigroup BTE (Eq. \eqref{eq:bte}),
\item the  multigroup LOQD equations  (Eq. \eqref{eq:mg-loqd}),
\item the  effective grey  LOQD equations  (Eq. \eqref{eq:gr-loqd}),
\item the MEB equation in effective grey form (Eq. \eqref{grey_meb}).
\end{itemize}

We discretized in time the high-order BTE, the hierarchy of the moment equations, and MEB  with a fully implicit temporal scheme over the given temporal mesh using  the backward Euler (BE) time-integration method.
Thus, we apply a one-stage time-step method to the multilevel  system of equations.
In a standard approach, this  system of equations is solved on each $n^{th}$ time step,
denoted as $\theta^n=(t^{n-1},t^{n}]$, iteratively to compute the solution
at  $t^{n}$ and then move to the next layer on time at $t^{n+1}$ \cite{PASE-1986,dya-aristova-vya-mm1996,aristova-vya-avk-m&c1999,dya-jcp-2019,pg-dya-jcp-2020,dya-jcp-2021}.
In this paper, we present a new iteration method defined for agglomerated sets of time steps from the original target grid in time.
These sets of time steps form coarse time intervals with a temporal subgrid defined by the original grid in time.
The multilevel system of equations approximated with the BE scheme over time steps $\theta^n$ is solved iteratively over the coarse intervals on the imbedded time subgid.
This iterative method involves (i) cycles of solving the high-order BTE over the set of time steps included in a coarse time interval and then
(ii) cycles of solving the hierarchy of low-order equations coupled with MEB equations over the same temporal subgrid.
This numerical method can be interpreted as an implicit multi-stage time integration method.
The developed iteration method is stable.
As the number of time steps included in the coarse time interval increases, the rate
of convergence decreases. However, the outer  iteration cycles over coarse time intervals still converge rapidly.
The  characteristic rate of convergence doesn't exceed 0.2.
The analysis of this method  enables us to study the effects of decoupling the high-order and low-order parts of multiphysics system of equations over  multiple time steps.
The elements of this new method have a potential for developing algorithm for parallel computations.

The remainder of this paper is organized as follows.
In Sec. \ref{sec:iter-meth} the multilevel iteration method for TRT problems is described.
Numerical results are presented in In Sec. \ref{sec:num-res}.
We conclude with a brief discussion in Sec.  \ref{sec:conc}.

\section{Formulation of the Iteration Method \label{sec:iter-meth}}
Consider that the TRT problem is discretized on a temporal mesh defined by $N$ time steps $\{\theta^n\}_{n=1}^N$ over the interval $t \in [0,t^\text{end}]$ such that
\begin{equation}
	\Big\{ \theta^n=(t^{n-1},t^n]\ |\ n=1,\ldots,N,\ t^0=0,\ t^N=t^\text{end} \Big\}.
\end{equation}

The standard strategy is to iterate over the entire hierarchy of equations \eqref{eq:bte} \& \eqref{eq:mg-loqd} - \eqref{eq:grip} at each time step $\theta^n$ separately.
The iterative scheme presented here instead performs iterations over collections of individual time steps together.
Let us define a subset of the discrete instants of time $\{\mathfrak{T}_b\}_{b=0}^B\subset\{t^n\}_{n=0}^N$:
\begin{equation}
	\Big\{ \mathfrak{T}_b=t^{N_b}\ |\ b=0,\dots,B,\quad 0=N_0<N_1<\dots<N_b=N \Big\}
\end{equation}
and split the temporal domain into $B$ {\it time blocks} $\{\Theta_b\}_{b=1}^B$ defined by
\begin{equation}
 	\Big\{ \Theta_b =  (\mathfrak{T}_{b-1},\mathfrak{T}_b]\ |\ b=1,\ldots,B \Big\}.
\end{equation}
Each time block is an amalgamated sequence of consecutive time steps $\Theta_b = \bigcup_{n=N_{b-1}+1}^{N_b} \theta^n$.
The number of time steps embedded in $\Theta_b$ is denoted $\mathfrak{N}_b=N_{b}-N_{b-1}$
The iterative algorithm presented here performs cycles over time blocks $\Theta_b$, effectively treating each time block
as a coarse time step.
The coupled solution of the BTE and LOQD equations is thus iteratively converged on time blocks instead of on each individual time step.
On each outer iteration,
the high-order BTE (Eq. \eqref{eq:bte}) is solved over all  time steps $\theta^n \in \Theta_b$
using the  temperatures $\{ T(t^n)\ |\ n =N_{b-1}+1,\ldots,N_b  \}$  over this collection $\Theta_b$ estimated on the previous iteration cycle.
The Eddington tensor is calculated for each time step involved in this process and stored until completion of the entire time block.
Afterwards the collected $\{\boldsymbol{\mathfrak{f}}_g(t^n)\ |\ n=N_{b-1}+1,\dots,N_b\}$ is used to close the LOQD equations for all time steps within the block $\Theta_b$.
The LOQD  and MEB Eqs. \eqref{eq:mg-loqd}-\eqref{grey_meb} are then solved over $\theta^n \in  \Theta_b$.
This yields the updated material temperature field $\{ T(t^n)\ |\ n =N_{b-1}+1,\ldots,N_b  \}$.
This iterative process is derived by noting that the BTE can be solved for all time steps in a given block provided the required initial condition $I_g|_{t=\mathfrak{T}_{b-1}}$ and the
iterative estimate of temperature field $\{ T(t^n)\ |\ n =N_{b-1}+1,\ldots,N_b  \}$.
Similarly, the Eqs. \eqref{eq:mg-loqd}-\eqref{grey_meb} can be solved for all time steps in the same block provided initial conditions for $E_g,\boldsymbol{F}_g,T$ at $t=\mathfrak{T}_{b-1}$ and the iterative estimate of  Eddington tensor  $\{ \boldsymbol{\mathfrak{f}}_g(t^n) \ | \ n =N_{b-1}+1,\ldots,N_b \}$.

The detailed iterative scheme is presented in Algorithm \ref{alg:serial}
where $j$ is the outer-iteration index.
The zeroth outer iteration involves solving the low-order system over all
time steps in $\Theta_b$
using the initial guess for the Eddington tensor
$\boldsymbol{\mathfrak{f}}_g= \frac{1}{3} \mathbb{I} $,
where $\mathbb{I}$ is the identity matrix.
During the iteration process, the multilevel system of low-order and MEB equations
is solved by means of two nested iteration cycles \cite{dya-jcp-2019}.
Convergence criteria are based on the function $\xi^{(j)}(E,t^n)=\| E^{(j)}(t^n) - E^{(j-1)}(t^n) \|_2$. We also denote $\bar{E}^{(j)}(t^n)=\| E^{(j)}(t^n) \|_2$. The specific norm for convergence is defined $\|\xi^{(j)}(E)\|_\infty=\max_{n\in[N_{b-1}+1  ,N_b]}\xi^{(j)}(E,t^n)$, and $\|\bar{E}^{(j)}\|_\infty=\max_{n\in[N_{b-1}+1  ,N_b]}\bar{E}^{(j)}(t^n)$. The same notations are used for norms of $T$.

\begin{algorithm}[ht!]
	{\footnotesize
		\SetAlgoLined
		\For{$b=1,\dots,B$}{
			$j=0$\\
			\For{$n=N_{b-1}+1, \dots,N_b$}{
				Given $\boldsymbol{\mathfrak{f}}^{(j)}_g(\! t^n\! )=\frac{1}{3}$: Solve LOQD \& MEB Eqs. \eqref{eq:mg-loqd} - \eqref{grey_meb} over $\theta^n$ for $ E_g^{(j)}(\! t^n\! ), \! \boldsymbol{F}_g^{(j)}(\! t^n\! ), \!  E^{(j)}(\! t^n\! ), \! \boldsymbol{F}^{(j)}(\! t^n\! ),  T^{(j)}(\! t^n\! )$
			}
			\While{ $\|\xi^{(j)}(E)\|_\infty > \epsilon\|\bar{E}^{(j)}\|_\infty,\ \ \|\xi^{(j)}(T)\|_\infty > \epsilon\|\bar{T}^{(j)}\|_\infty$ }{
				$j=j+1$\\ \vspace{0.2cm}
				$\bullet$ \underline{BTE  for $t \in \Theta_b$} \\
				Given $\{ T^{(j-1)}(t^n)\ |\ n\in[N_{b-1}+1  ,N_b] \}$   \&  $I_g^{(j)}(t^{N_{b-1}})$:   \\ \vspace{0.1cm}
				\For{$n=N_{b-1}+1,\dots,N_b$}{
					Given ${T}^{(j-1)}(t^n)$: Solve BTE \eqref{eq:bte}  over $\theta^n$ for $I_g^{(j)}(t^n)$\\ \vspace{0.1cm}
					Calculate $\boldsymbol{\mathfrak{f}}^{(j)}_g(t^n)$ with Eq. \eqref{eq:et}
				}
				\vspace{0.2cm}
				$\bullet$ \underline{LOQD \& MEB eqs.  for $t \in \Theta_b$} \\
				Given $\{ \boldsymbol{\mathfrak{f}}_g^{(j)}(t^n)\ |\ n\in[N_{b-1}+1,N_b] , g=1,\ldots,G\}$, $E_g^{(j)}(t^{N_{b-1}})$, $\boldsymbol{F}_g^{(j)}(t^{N_{b-1}})$, $T^{(j)}(t^{N_{b-1}})$: \\
				$s=0$\\
				\For{$n=N_{b-1}+1,\dots,N_b$}{
					\While{$\| E^{(s,j)}(t^n) - E^{(s-1,j)}(t^n) \|_2 > \tilde \epsilon\|E^{(s,j)}(t^n)\|_2,\ \ \| T^{(s,j)}(t^n) - T^{(s-1,j)}(t^n) \|_2 > \tilde \epsilon\|T^{(s,j)}(t^n)\|_2$ }{
						$s=s+1$\\
						Given $\boldsymbol{\mathfrak{f}}^{(j)}_g(t^n)$ and $T^{(s-1,j)}$: Solve multigroup LOQD Eqs. \eqref{eq:mg-loqd} over $\theta^n$ for
						$E_g^{(s,j)}(t^n),\boldsymbol{F}_g^{(s,j)}(t^n)$\\
						Compute grey opacities and grey Eddington tensor (Eqs. \eqref{eq:grip2} \& \eqref{eq:grip})\\
						Solve grey LOQD eqs. \eqref{eq:gr-loqd} coupled with MEB eq. \eqref{grey_meb} over $\theta^n$  to compute $E^{(s,j)}(t^n),\boldsymbol{F}^{(s,j)}(t^n),  T^{(s,j)}(t^n)$
						
					}
					
			}	}
		}
	}
	\caption{
		\label{alg:serial}
		Iterative method with outer iteration cycles
		of  high-order and low-order equations   over multiple time steps  }
\end{algorithm}

This iterative algorithm can be interpreted as one which solves a multi-step implicit temporal discretization of the TRT problem, similar to diagonally-implicit Runge-Kutta methods.
Consider the BTE discretized in time with the backward-Euler scheme on the grid of time steps $\{\theta^n\}_{n=1}^N$
\begin{equation} \label{eq:bte-be}
	I^n_g = \tau^nH^n_g(I^n_g,T^n) + I^{n-1}_g,
\end{equation}
where $\tau^n=c\Delta t^n$, $\Delta t^n=t^n-t^{n-1}$ and
\begin{equation}
	H^n_g(I^n,T^n) = \varkappa_g(T^n)B_g(T^n) - \boldsymbol{\Omega}\cdot\boldsymbol{\nabla}I_g^n - \varkappa_g(T^n)I_g^n.
\end{equation}
Taking Eq. \eqref{eq:bte-be} at the last time step in a given time block $(n=N_b)$, and recursively substituting Eq. \eqref{eq:bte-be} for $\{I_g^{n-1}\}_{n=N_{b-1}+1}^{N_b}$ yields the following:
\begin{equation} \label{eq:bte-be-sum}
	I^{N_b}_g = I^{N_{b-1}}_g + \sum_{n=N_{b-1}+1}^{N_b}\tau^nH^n_g(I^n_g,T^n).
\end{equation}
Let us define the coarse time-block-step $\Delta\mathfrak{T}_b = \mathfrak{T}_b-\mathfrak{T}_{b-1}$.
Eq. \eqref{eq:bte-be-sum} can be rewritten in the form of a multi-step method with step length $\Delta\mathfrak{T}_b$ as
\begin{equation}
	I^{N_b}_g = I^{N_{b-1}}_g + \Delta\mathfrak{T}_b\sum_{n=N_{b-1}+1}^{N_b}\frac{\tau^n}{\Delta\mathfrak{T}_b}H^n_g(I^n_g,T^n).
\end{equation}
The same process can be used to derive this multi-step temporal discretization over time blocks for the LOQD and MEB equations.

\section{Numerical Results \label{sec:num-res}}
The iterative scheme is analyzed with numerical testing on the classical Fleck-Cummings (F-C) test \cite{fleck-1971} in 2D Cartesian geometry. Standard formulation of the test is used.
The test domain is a $6\times 6$ cm homogeneous slab with spectral opacity $\varkappa_\nu = \frac{27}{\nu^3}(1-e^{-\nu/T})$.
The left boundary is subject to incoming radiation at a temperature of $T^\text{in}=1$ KeV at black-body spectrum. All other boundaries are vacuum.
The initial temperature of the slab is $T^0=1$ eV. The material energy density of the slab is a linear one $\varepsilon=c_vT$ where $c_v=0.5917 a_R T^\text{in}$ and $a_R$ is Stefan's constant.
We consider a time interval of $t\in[0,6 \, \text{ns}]$, discretized into 300 uniform time steps $\Delta t = 2\times 10^{-2}$ ns. The phase space is discretized using a $20\times 20$ uniform orthogonal spatial grid, 17 frequency groups \cite{pg-dya-jcp-2020}, and 144 discrete directions.
The Abu-Shumays angular quadrature set q461214 is used \cite{abu-shumays-2001}.
The implicit backward-Euler time integration scheme is used to discretize all equations in time. The BTE is discretized in space with the method of long characteristics,  and all low-order equations use a second-order finite-volumes scheme \cite{jcm-dya-mc2023}.

We consider cases with time blocks of lengths $\Delta\mathfrak{T}_b=0.02,\ 0.04,\ 0.1,\ 0.2,\ 0.5,\ 1,\ 2,\ 3,\ 6$ ns
(i.e. $\mathfrak{N}_b = 1, 2, 5, 10, 25, 50, 100, 150, 300$, respectively).
Note that when
$\Delta\mathfrak{T}_b=0.02$ ns, each time block is simply one time step ($\mathfrak{N}_b=1$).
Thus, this is the case of standard iteration scheme with outer iteration cycle
over each given time step.
In the case $\Delta\mathfrak{T}_b=6$ ns,
on each outer iteration cycle  the high-order BTE and low-order equations are solved
over the whole time range of the problem, i.e. $\Theta_1=[0,t^\text{end}]$.
The F-C test is solved using the presented iterative scheme and compared against the standard iteration scheme \cite{dya-jcp-2019}.

\begin{figure}[ht!]
	\centering
	\includegraphics[width=0.65\textwidth]{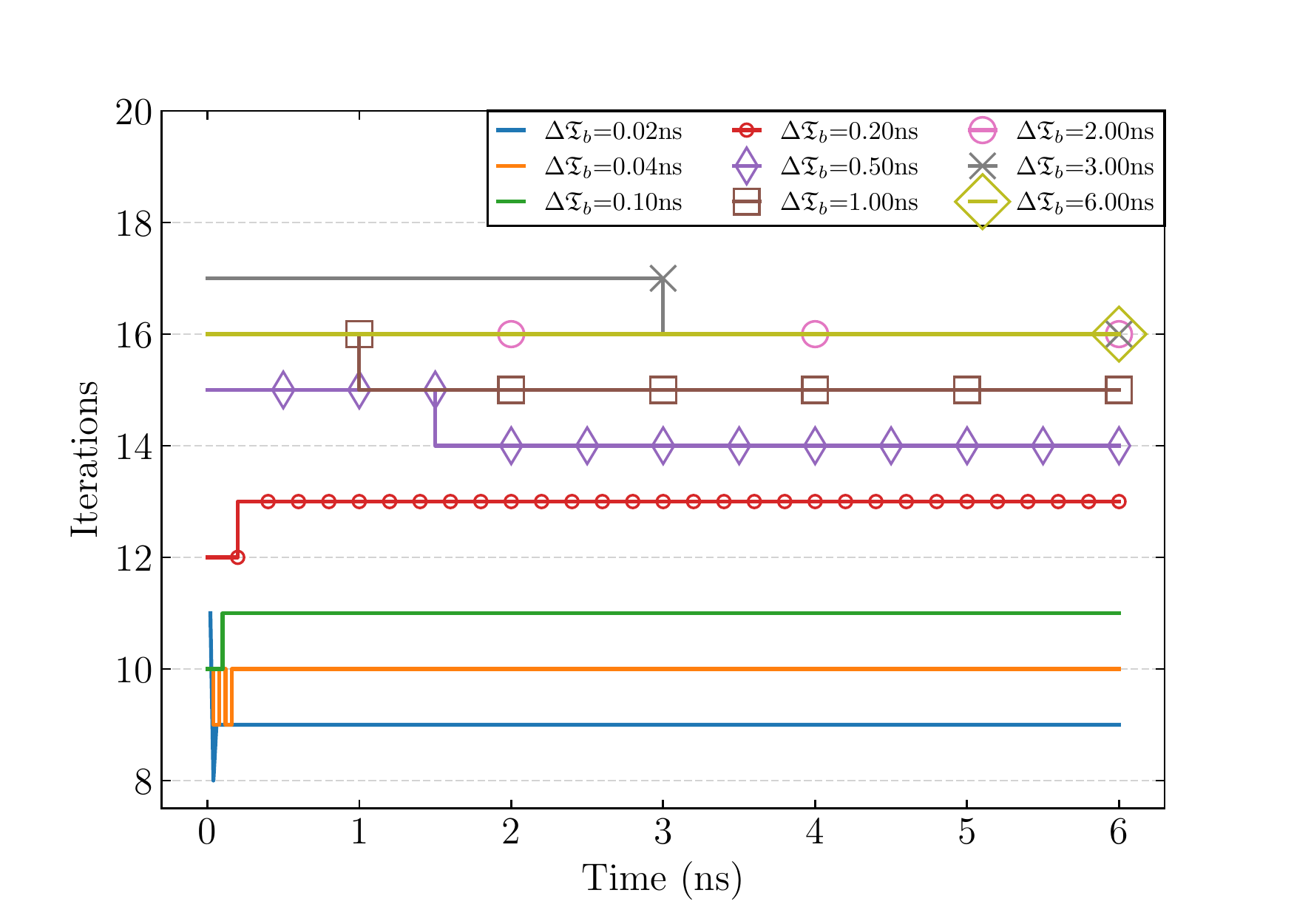}
	\caption{Iterations $(j)$ taken at each time block
		to reach a relative convergence level of $\epsilon=10^{-14}$. Blocks designated by markers for sufficiently long block times. \label{fig:its}
	}
\end{figure}

We use a very tight convergence criteria of $\epsilon=10^{-14}$.
This enables us to   analyze  the convergence behavior of the proposed iterative scheme   over
the large range of iterative errors up to very small orders of magnitude
and  avoid missing iterative stagnation and noise effects.
Figure \ref{fig:its} plots the iteration counts per block $b$ to reach the required convergence level.
For sufficiently long time blocks, markers are placed at the end of each block interval.
As time blocks become larger, the iteration counts tend to increase.
The iteration counts for $\Delta\mathfrak{T}_b=0.02$ ns are the same as for the standard iteration scheme, and therefore the presented scheme requires more iterations (per block) than the standard scheme for all $\Delta\mathfrak{T}_b>0.02$ ns.
These effects stem from the fact that the high-order and low-order problems are coupled over entire time-block cycles, instead of over each time step.

\begin{figure}[h!]
	\centering
	\subfloat[$\Delta\mathfrak{T}_b=0.02\ \text{ns},\ \mathfrak{N}_b=1$]{\includegraphics[width=0.47\textwidth]{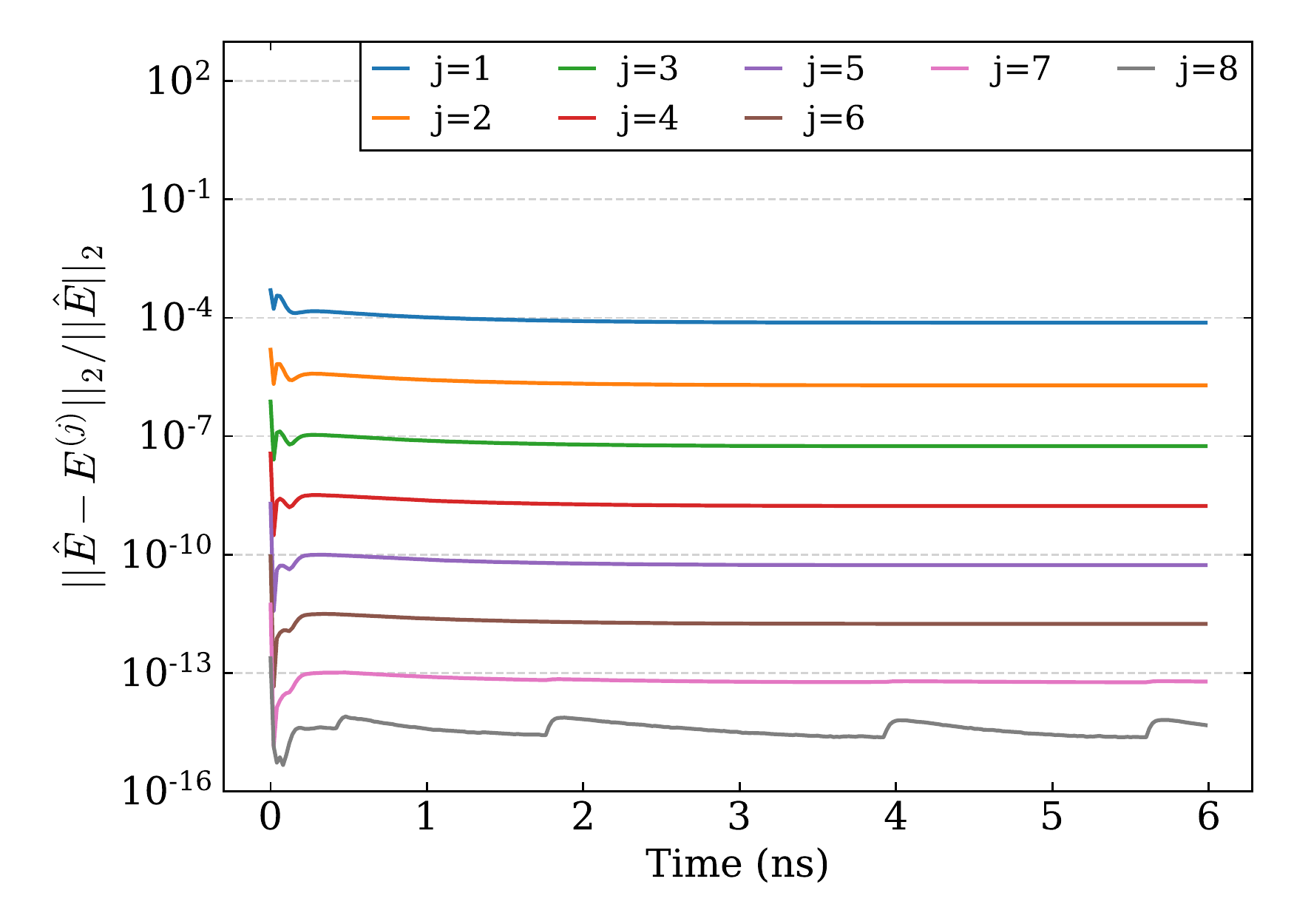}}
	\subfloat[$\Delta\mathfrak{T}_b=0.10\ \text{ns},\ \mathfrak{N}_b=5$]{\includegraphics[width=0.47\textwidth]{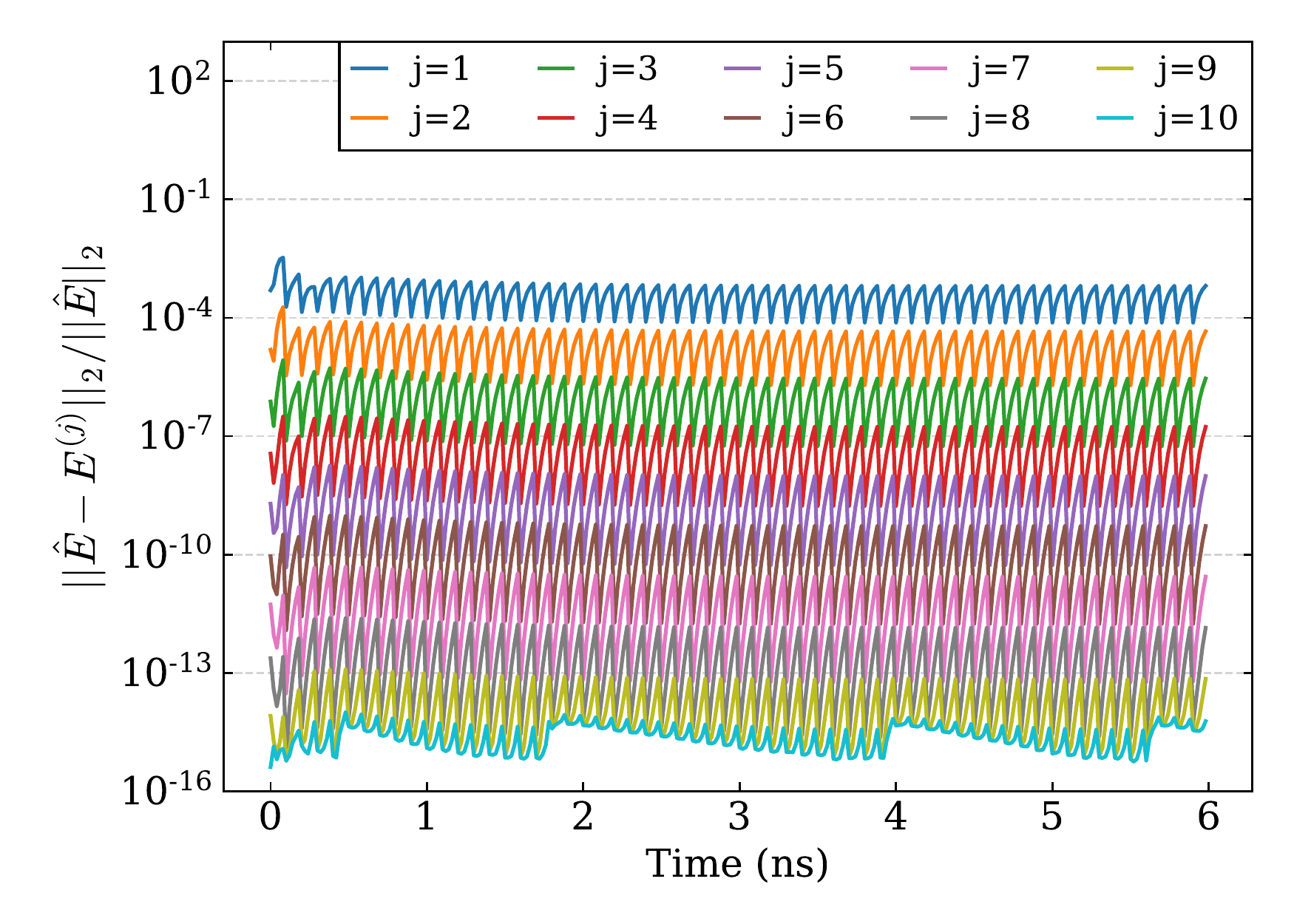}}\\[5pt]
	\subfloat[$\Delta\mathfrak{T}_b=0.5\ \text{ns},\ \mathfrak{N}_b=25$]{\includegraphics[width=0.47\textwidth]{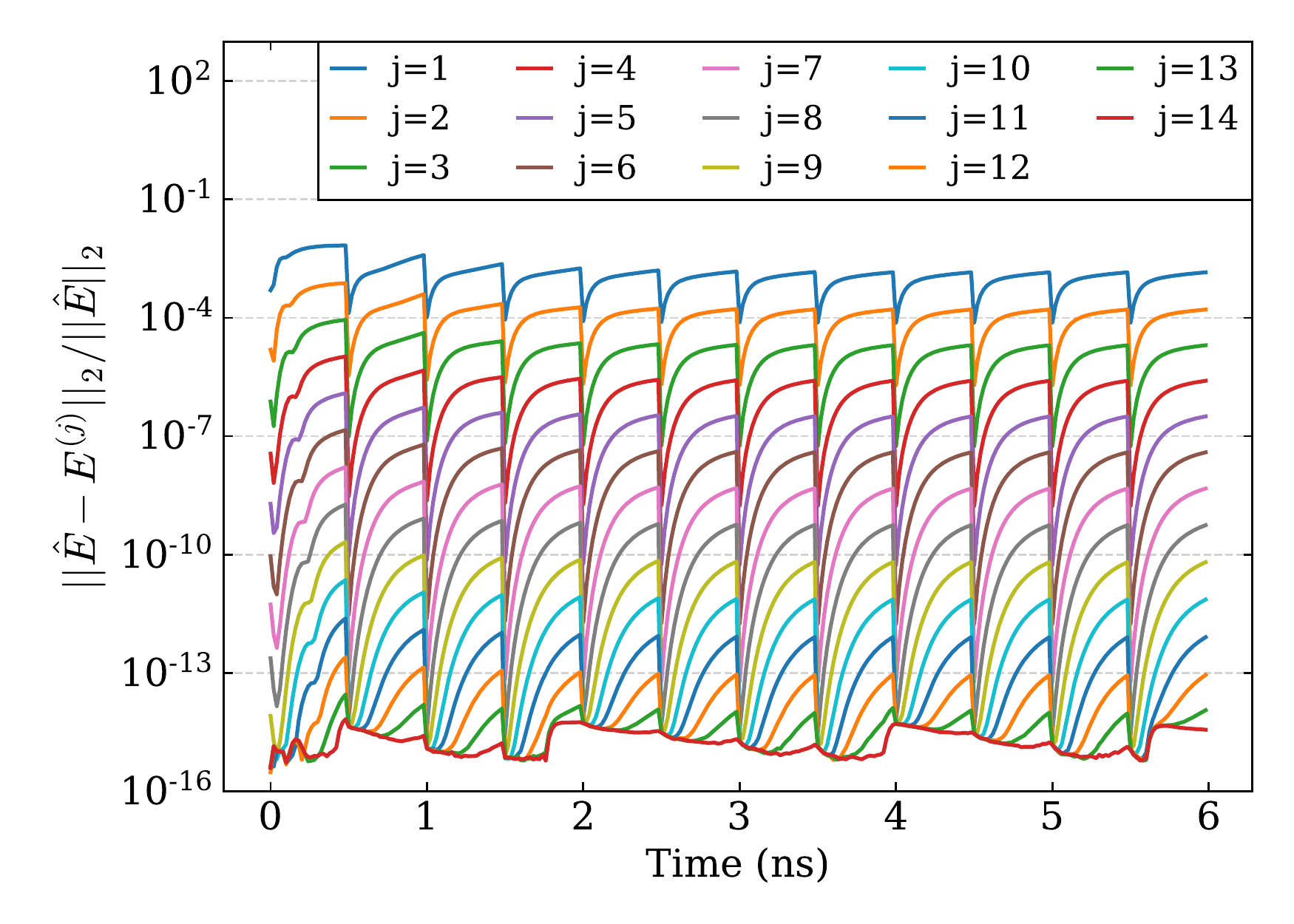}}
	\subfloat[$\Delta\mathfrak{T}_b=1\ \text{ns},\ \mathfrak{N}_b=50$]{\includegraphics[width=0.47\textwidth]{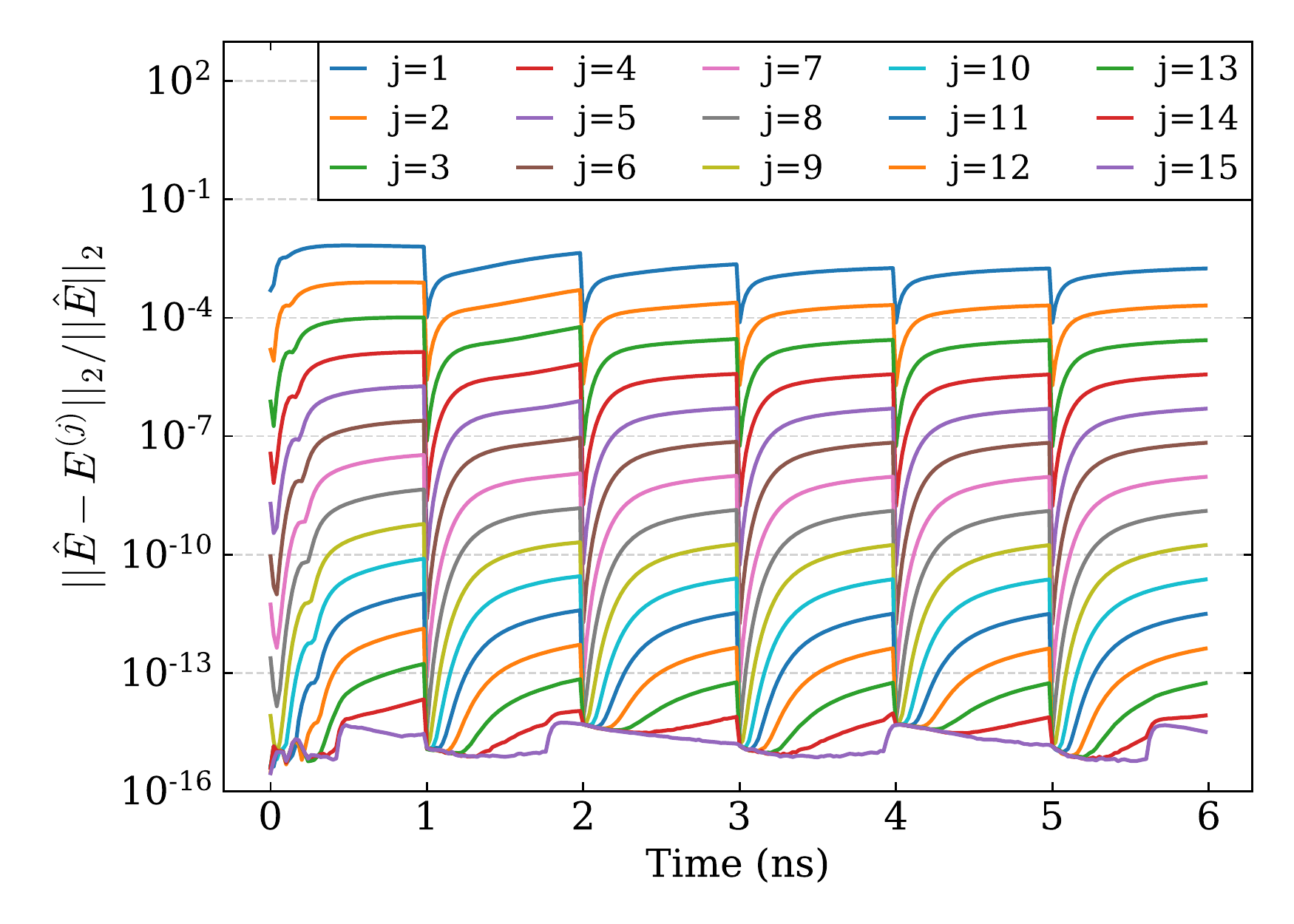}}\\[5pt]
	\subfloat[$\Delta\mathfrak{T}_b=2\ \text{ns},\ \mathfrak{N}_b=100$]{\includegraphics[width=0.47\textwidth]{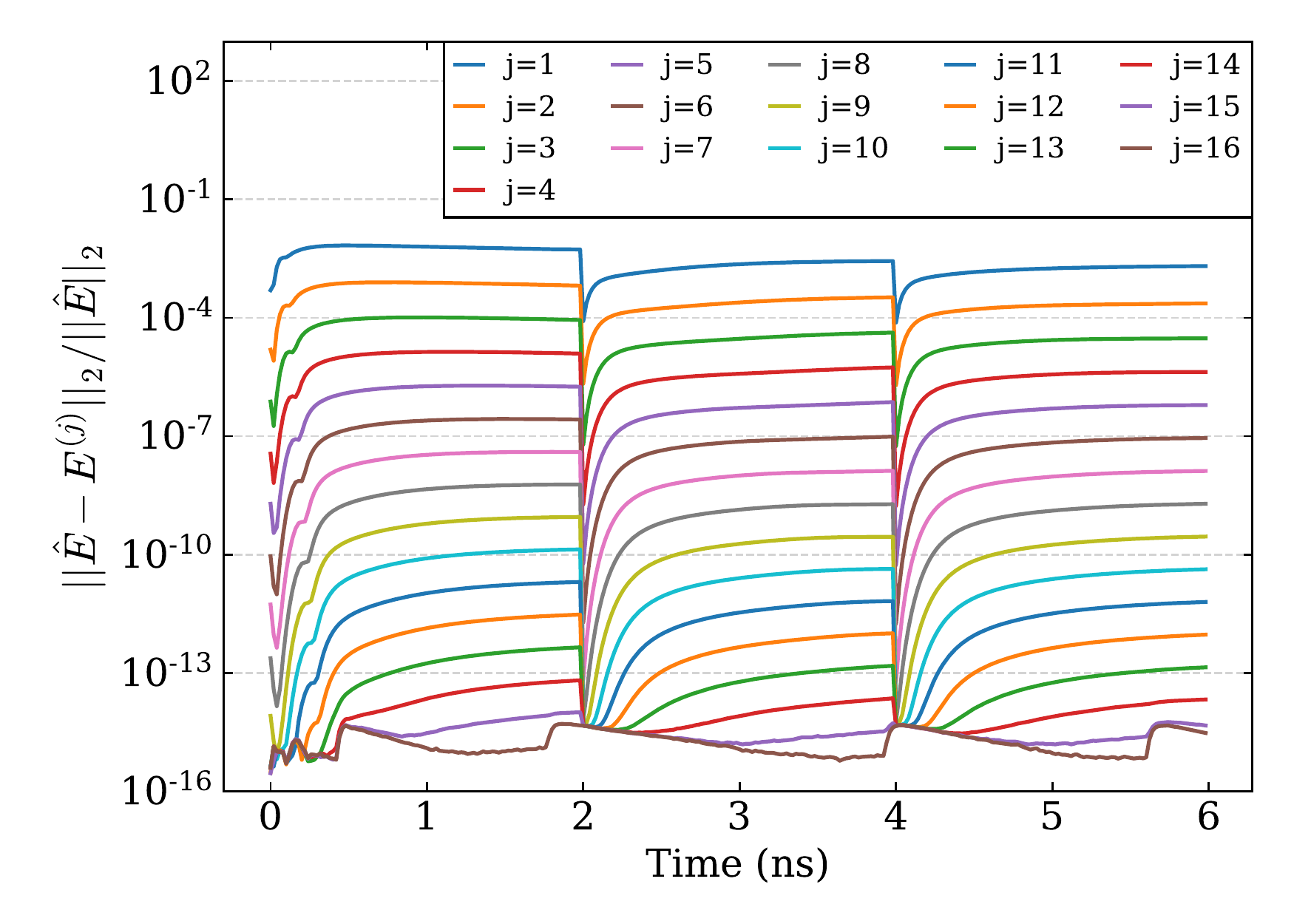}}
	\subfloat[$\Delta\mathfrak{T}_b=3\ \text{ns},\ \mathfrak{N}_b=150$]{\includegraphics[width=0.47\textwidth]{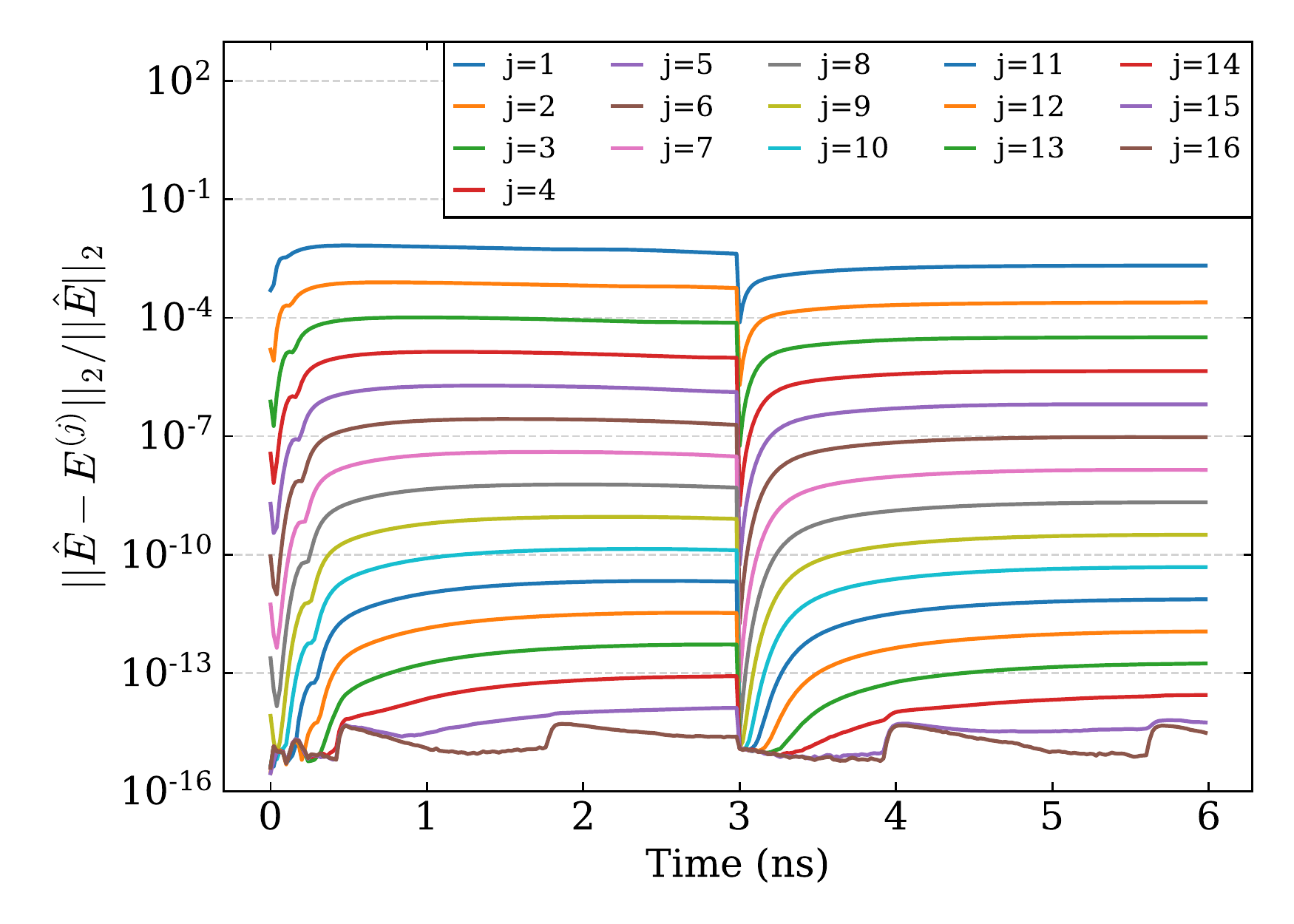}}
	\caption{Relative error in $E$
		computed with several $\Delta\mathfrak{T}_b$ w.r.t. the reference solution.
		\label{fig:conv}
	}
\end{figure}

Figure \ref{fig:conv} plots the iterative error in the total radiation energy density $E$ at each time step and iteration for several cases, including $\Delta\mathfrak{T}_b=0.02,\ 0.1,\ 0.5,\ 1,\ 2,\ 3$ ns. The standard case with $\mathfrak{N}_b=1$ is included for reference.
Note that errors in the solution at each iteration $(j)$ are calculated with respect to the
reference  solution computed by means of the standard scheme, (i.e. with 
$\Delta\mathfrak{T}_b =0.02$ns)  denoted as $\hat{E}$.
The developed iterative scheme's solution is shown to converge
to the reference solution.
Note that similar results are obtained for errors in $T$, and all other tested values for $\Delta\mathfrak{T}_b$. Errors are seen to converge uniformly at each iteration. The errors in each block sharply increase over the first several time steps, then level off.

\begin{figure}[h]
	\centering
	\subfloat[$\Delta\mathfrak{T}_b=0.04\ \text{ns},\ \mathfrak{N}_b=2$]{\includegraphics[width=0.5\textwidth]{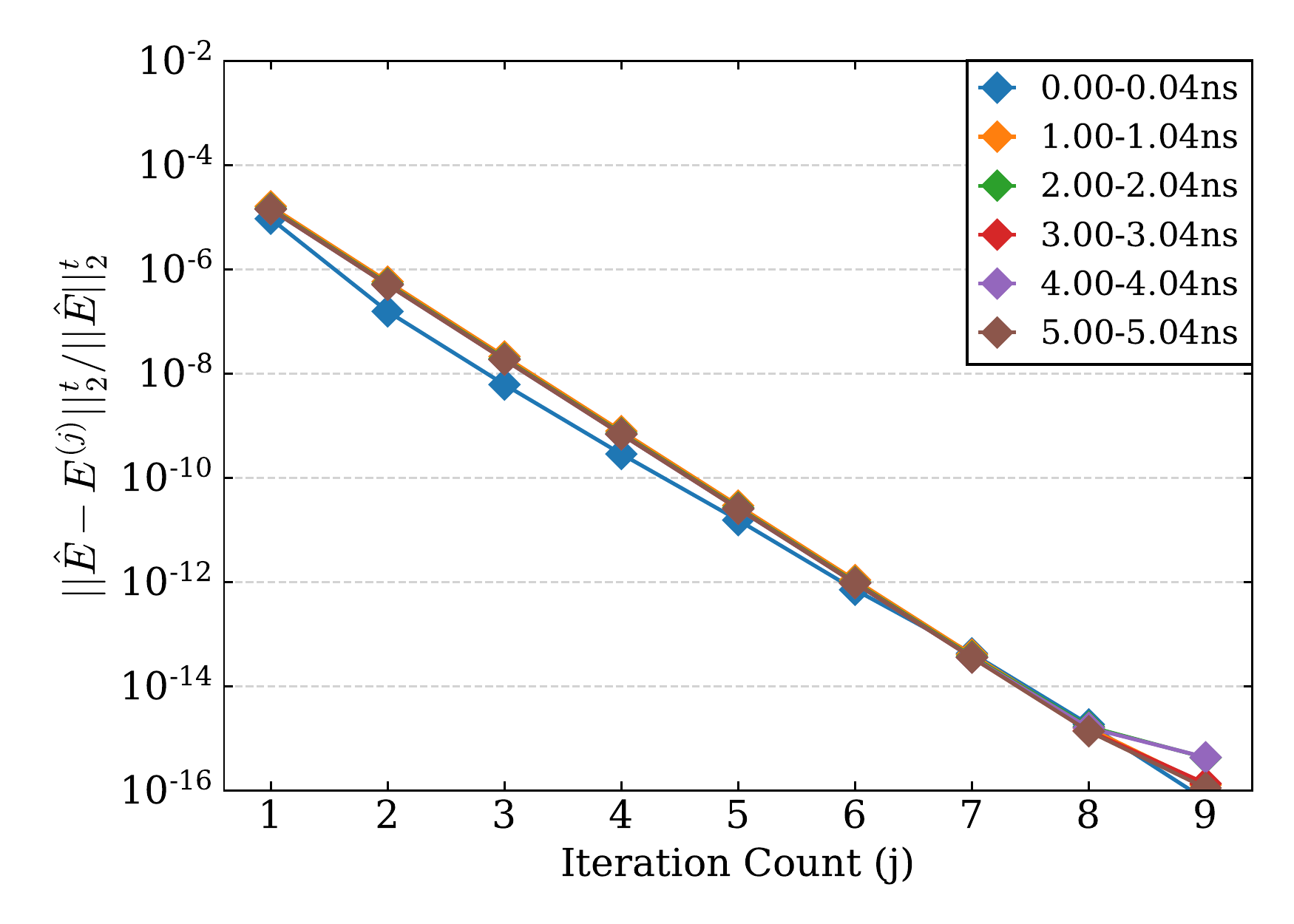}}
	\subfloat[$\Delta\mathfrak{T}_b=0.1\ \text{ns},\ \mathfrak{N}_b=5$]{\includegraphics[width=0.5\textwidth]{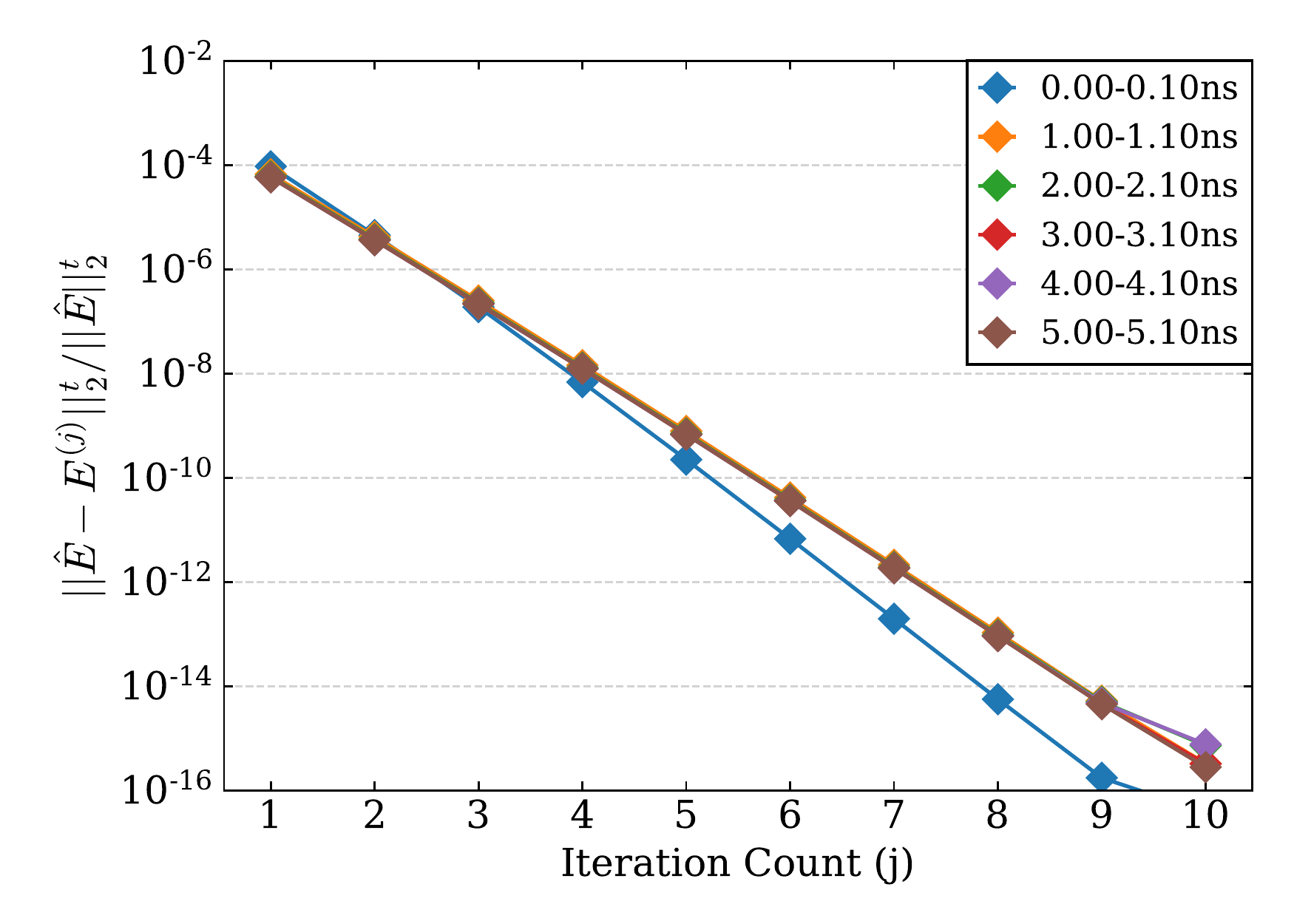}}\\
	\subfloat[$\Delta\mathfrak{T}_b=0.5\ \text{ns},\ \mathfrak{N}_b=25$]{\includegraphics[width=0.5\textwidth]{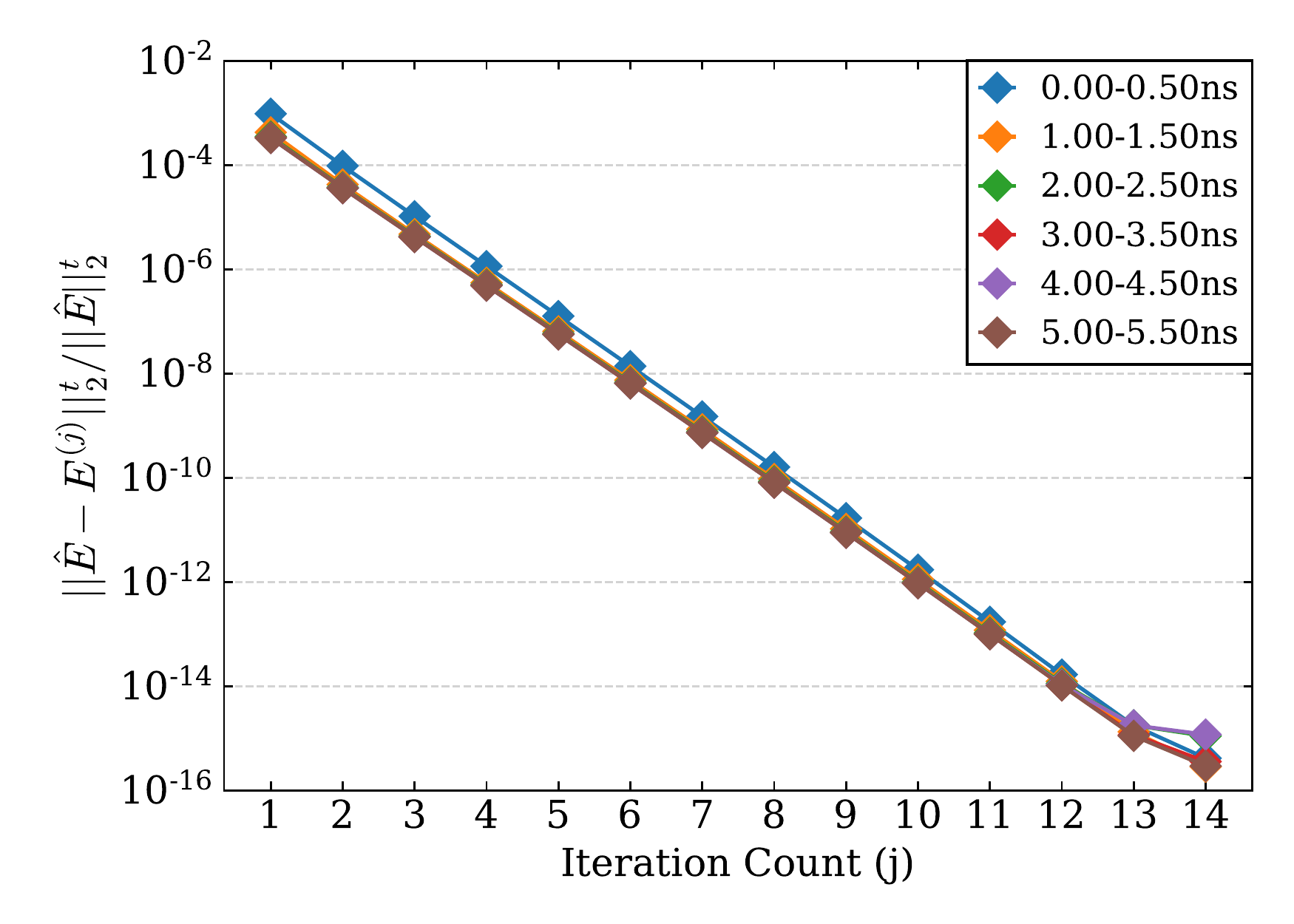}}
	\subfloat[$\Delta\mathfrak{T}_b=1\ \text{ns}\ \mathfrak{N}_b=50$]{\includegraphics[width=0.5\textwidth]{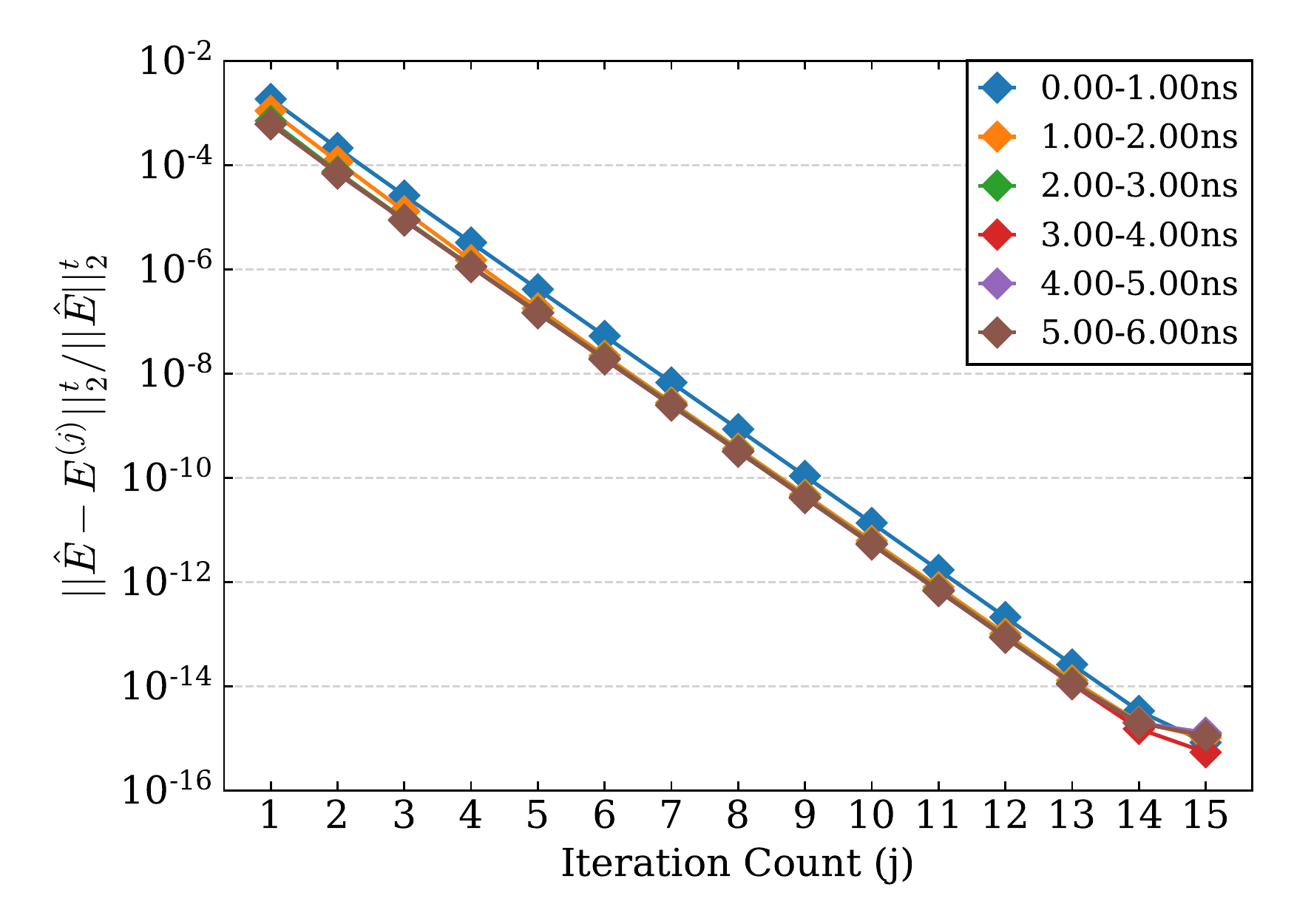}}
	\caption{Convergence plots of relative error in $E$ w.r.t. the reference solution vs iteration for several select time intervals and $\Delta\mathfrak{T}_b$ values.
		\label{fig:spr}
	}
\end{figure}

Figure \ref{fig:spr} plots the iterative error in $E$ using $\Delta\mathfrak{T}_b=0.04,\ 0.1,\ 0.5,\ 1$ ns w.r.t. the reference solution computed by the standard iterative scheme. The plot is formatted so that results are graphed against the iteration count, and each line corresponds to a different time block. Errors are calculated in the norm $\|\cdot\|_2^t$, which is the 2-norm over space and time for the temporal interval of a given time block. The blocks chosen on the graphs are representative of the overall behavior.
Table~\ref{tab:spr} displays the average rate of convergence of iterations for both $E$ and $T$ estimated by averaging over all time intervals for each tested $\Delta\mathfrak{T}_b$. The estimated average rate of convergence for the $j^\text{th}$ iteration of a given time interval is calculated as: $\rho^{(j)}_E={\|\hat{E}-E^{(j)}\|_2^t}/{\|\hat{E}-E^{(j+1)}\|_2^t}$. The values shown in Table \ref{tab:spr} have been averaged over all iterations and time blocks. The average rate of convergence is observed to increase with $\Delta\mathfrak{T}_b$ and eventually level off at around $\Delta\mathfrak{T}_b=1.00$ ns.

\begin{table}[h]
	\centering
	\vspace{-.35cm}
	\caption{Estimated average rate of convergence of iterations for both $E$ and $T$ averaged over all time intervals \label{tab:spr}}
\vspace{0.2cm}
	\begin{tabular}{|c|c|c|c|}
		\hline
		$\Delta\mathfrak{T}_b$ (ns) & $\mathfrak{N}_b$ & $\rho_E$ & $\rho_T$\\\hline
		0.02 & 1 & 0.042 & 0.035\\\hline
		0.04 & 2 & 0.068 & 0.049\\\hline
		0.10 & 5 & 0.067 & 0.058\\\hline
		0.20 & 10 & 0.128 & 0.100\\\hline
		0.50 & 25 & 0.158 & 0.136\\\hline
		1.00 & 50 & 0.171 & 0.154\\\hline
		2.00 & 100 & 0.194 & 0.178\\\hline
		3.00 & 150 & 0.177 & 0.167\\\hline
		6.00 & 300 & 0.159 & 0.156\\\hline
	\end{tabular}
\end{table}

\section{Conclusions \label{sec:conc}}
In this paper, a new projective-iterative scheme for TRT problems is presented.
The algorithm performs iteration cycles which solve (i) high-order equations and (ii) low-order equations
over collections of multiple time steps. This is achieved by decoupling the high- and low- order systems in time.
The numerical results show this iteration method to be stable
in solving the classical TRT problem with the radiation Marshak wave for a very  large range of lengths of   time blocks.
The method converges even with a cycle over a single time block that covers the whole time interval of the problem.
As the length of time blocks is increased, the number of outer iterations (cycles) increases.

A potential advantage of this iterative scheme is the possibility for parallelization in time.
Since the high- and low-order problems are solved separately from one another during the iterations of a given time block, they could be solved in parallel to one another.
Further development and analysis of the methodology is required to fully investigate this possibility. Some of the primary questions to answer include how often to communicate between high- and low- order equations the information they require, and optimal sizing of the coarse time blocks.
The sizing of time blocks will affect both the amount of memory required to store data between communications, how often communication must be performed, and how much computational load can be distributed among the different processes.

\section*{ACKNOWLEDGEMENTS}
 Los Alamos Report LA-UR-22-31465.
The  work of the first author (JMC) was supported by the U.S. Department of Energy through the Los Alamos National Laboratory. Los Alamos National Laboratory is operated by Triad National Security, LLC, for the National Nuclear Security Administration of U.S. Department of Energy (Contract No. 89233218CNA000001).

\bibliographystyle{elsarticle-num}
\bibliography{coale-anistratov-iter-scheme-cycles-MTS-TRT-preprint-arXiv}

\end{document}